\documentclass{amsart}

\usepackage{amssymb}

\usepackage[english,francais]{babel}

\newtheorem{lemme}{Lemme}
\newcommand{\parag}[1]{\par\smallskip\noindent\textbf{{#1}.} }
\def \CM{\mathbb{C}}
\def \exp {{\rm exp\,}}

\def \id {{\rm id\,}}
\def\dt{\delta}

\def\e{\varepsilon}  
\def\g{\gamma}

\def \s{\sigma}

\def \to{\longrightarrow}

\def \alg{\mathfrak{g}}

\def \< {{\langle }}
\def \> {{\rangle }}

\begin{document}
\title[Actions de groupes]{Un th{\'e}or{\`e}me sur les actions\\
  de groupes de dimension infinie}

\author{Jacques F{\'e}joz}
\email{fejoz@math.jussieu.fr}
\thanks{Le premier auteur a b{\'e}n{\'e}fici{\'e} d'un financement de l'ANR
  (Project BLANC07-3\_187245, Hamilton-Jacobi and Weak KAM Theory).}
\author{Mauricio Garay}
\email{garay91@gmail.com}
\maketitle

\begin{abstract}
  \selectlanguage{francais} 

  L'objet de cette note est de donner, dans un cadre analytique, un
  crit{\`e}re infinit{\'e}simal pour qu'un espace vectoriel soit localement
  homog{\`e}ne sous l'action d'un groupe. Notre approche diff{\`e}re de celles
  qui recourent {\`a} un th{\'e}or{\`e}me d'inversion locale
  (e.g. \cite{Mos66a,Zeh76b} ou \cite[Th{\'e}or{\`e}me~4.2.5]{Ser72}), par
  l'usage crucial de la structure de groupe sous-jacente. Ceci permet,
  en particulier, de remplacer l'estimation pour l'inverse de l'action
  de l'alg{\`e}bre de Lie en un plan tangent arbitraire, par une
  estimation sur les vecteurs tangents en l'origine. Notre
  d{\'e}monstration est bas{\'e}e sur la m{\'e}thode it{\'e}rative utilis{\'e}e par
  Kolmogorov et Arnold dans leur d{\'e}monstration du th{\'e}or{\`e}me des tores
  invariants \cite{Kol54,Arn63}. Le th{\'e}or{\`e}me de cette note sera
  utilis{\'e} dans des travaux ult{\'e}\-rieurs.

  \medskip

  \selectlanguage{english}
  \noindent{\bf A theorem on infinite dimensional group actions. }

  \noindent\textsc{Abstract.}
  We give an infinitesimal criterion, in the analytic setting, for a
  vector space to be locally homogeneous under some group action.  Our
  approach differs from those which resort to an inverse function
  theorem (e.g. \cite{Mos66a,Ser72,Zeh76b} or
  \cite[Th{\'e}or{\`e}me~4.2.5]{Ser72}), because we use the underlying group
  structure in an essential way. In particular, this allows to replace
  the estimate of the inverse map of the Lie algebra action at an
  arbitrary tangent plane, by an estimate of the vectors tangent at
  the origin. Our proof relies on the iterative method used by
  Kolmogorov and Arnold in their proof of the invariant tori
  theorem~\cite{Kol54,Arn63}. The theorem of this note will be used in
  subsequent works.
\end{abstract}

\selectlanguage{francais}

\parag{1} Un \emph{espace vectoriel} $E$ est dit {\em {\'e}chelonn{\'e}} s'il
est r{\'e}union des {\'e}l{\'e}m{\'e}nts d'une famille d{\'e}croissante de sous-espaces de
Banach $(E_s)$, $s \in ]0,1[$, telle que les inclusions $E_{s+\sigma}
\subset E_s$, $\s>0$, soient de norme au plus~$ 1$.  Nous noterons
$B^{E}_s$ la boule unit{\'e} de $E_s$ centr{\'e}e en l'origine
et $B_E$ l'union de ces boules.\\
Une application lin{\'e}aire $u : E \to F$ entre deux espaces vectoriels
{\'e}chelonn{\'e}s est \emph{$k$-born{\'e}e} si~:
$$N(u):=\sup_{s< s+\sigma, \; x \in E_{s+\sigma}} \sigma^k
\frac{| u(x)|_s}{\phantom{u}|x|_{s+\sigma}} < +\infty.$$ 

\smallskip De m{\^e}me, une \emph{vari{\'e}t{\'e} topologique} $X$ est dite
\emph{{\'e}chelonn{\'e}e} si elle est munie d'une famille croissante de
vari{\'e}t{\'e}s diff{\'e}rentielles banachiques $X_s \subset X$ model{\'e}es sur les
{\'e}l{\'e}ments d'un espace vectoriel {\'e}chelonn{\'e} $E$. 

\parag{2} Supposons qu'un ensemble $G=\bigcup_s G_s$ soit muni {\`a} la
fois d'une structure de groupe et de vari{\'e}t{\'e} {\'e}chelonn{\'e}e . La vari{\'e}t{\'e}
$G_s$ n'{\'e}tant pas n{\'e}c{\'e}ssairement un groupe, les conditions de
compatibilit{\'e} sont plus subtiles que dans le cas des groupes de Lie.\\
Par analogie avec la dimension finie, notons
$$\exp : 2B^\alg \rightarrow G, \quad \xi
\mapsto e^\xi, \quad 0 \mapsto 1$$ un param{\'e}trage local de $G$ par un
hom{\'e}omorphisme commutant aux inclusions des {\'e}chelles. Nous notons
$\log$ l'inverse {\`a} gauche de cette application. Nous dirons que la
structure de vari{\'e}t{\'e} {\'e}chelonn{\'e}e $(G_s)$ fait de $G$ un \emph{groupe
  {\'e}chelonn{\'e}} s'il existe $\kappa>0$ tel que la loi interne induise des
applications~:
$$\quad
B^\alg_{s+\sigma} \times \sigma B^\alg_{s} \rightarrow
\exp\left((1+\sigma)B^\alg_{s}\right) \subset G_s, \quad 
(\xi,\eta) \mapsto e^\xi e^{\eta}$$  
v{\'e}rifiant les conditions~:
$$\left\{
  \begin{array}[c]{lll}
    |\log \left(e^\xi \, e^\eta \right)|_s &\leq &|\xi|_{s+\sigma} +
    |\eta|_s \\
    |\log e^\xi e^\eta - \xi - \eta|_s &\leq &\kappa \sigma^{-1}
    |\xi|_{s+2\sigma} |\eta|_{s} .
  \end{array} \right.$$
pour tous  $ \xi \in B^\alg_{s+2\sigma}, \, \eta \in
    \sigma B^\alg_{s}$.
\parag{3} Une \emph{action de groupe {\'e}chelonn{\'e}} $G=\bigcup G_s$
sur un espace vectoriel {\'e}chelonn{\'e} $E=\bigcup E_s$ est une action du
groupe $G$ sur $E$ qui induit des applications de classe $C^1$~:
$$\quad
(\sigma B^\alg_{s}) \times E_{s+\sigma} \rightarrow E_{s},
\quad (\xi,x) \mapsto e^\xi\,x.$$ 
Une telle action induit une action infinit{\'e}simale de $\alg$~:
$$\xi : E_{s+\sigma} \rightarrow E_{s}, \quad \xi\, x := 
\left.\frac{d}{dt}\right|_{t=0} \left( e^{t\xi}\, x \right)$$
pour tout $ \xi\in \alg_s$. Nous dirons que l'action
  v{\'e}rifie {\em une condition $(A_c)$},  $c\geq 1$, si~:
  $$\leqno{(A_c)}\quad
\left|\left(e^{\xi}\right)^{-1} \, \left( \xi \, 0_E \right)
\right|_{s} \leq c \, \sigma^{-1} |\xi|_{s+2\sigma}^2,{\rm ,\ pour\ tout\ }
  \xi\in   B^\alg_{s+2\sigma}$$
o{\`u} $0_E$ d{\'e}signe l'origine dans $E$.
\parag{4} 
L'exemple que nous avons en t{\^e}te est le suivant. Soit $U_s$ un syst{\`e}me
fondamental de voisinages d'un compact $K \subset \CM^n$.  Notons
$E_s$ l'espace des fonctions holomorphes sur $U_s$ et continues sur
l'adh{\'e}rence de $U_s$. Le groupe $\widehat G$ des transformations
biholomorphes qui laissent $K$ invariant peut {\^e}tre muni d'une
structure de groupe {\'e}chelonn{\'e} qui satisfait une condition $(A_c)$ et
qui d{\'e}finit une action {\'e}chelonn{\'e} sur $E=\bigcup E_s$.  La structure de
groupe {\'e}chelonn{\'e} est donn{\'e}e par la param{\'e}trisation $e^\xi := \id +
\xi$ et l'action par $gx = (a+x) \circ g^{-1} - a$ avec $g=e^\xi$.\\
Dans des travaux ult{\'e}rieurs, nous consid{\`e}rerons le cas o{\`u} $n=2d+l$,
$K$ est de la forme $L \times \{ 0\}$ o{\`u} $L \subset \CM^{2d}$ est la
partie r{\'e}elle d'une vari{\'e}t{\'e} lagrangienne complexe analytique, l'espace
$\CM^l$ est un espace de param{\`e}tres et $G \subset \widehat G$ est le
groupe des automorphismes de Poisson de $\CM^{2d} \times \CM^l$ au
voisinage de $K$.

\medskip\noindent\textsc{Th{\'e}or{\`e}me. }{\it
  Soit $G=\bigcup_s G_s$ un groupe {\'e}chelonn{\'e} agissant sur un espace
  vectoriel {\'e}chelonn{\'e} $E=\bigcup_s E_s$. Supposons d'une part que
  l'action v{\'e}rifie une condition $(A_c)$, et que d'autre part
  l'application lin{\'e}aire $\rho : \alg \to E,\ \xi \mapsto \xi \, 0_E$
  poss{\`e}de un inverse {\`a} droite $j$ qui soit $k$-born{\'e}. L'espace $E$ est
  alors localement $G$-homog{\`e}ne, plus pr{\'e}cis{\'e}ment~:
  $$\e \, B^E_{s} \subset \exp
  \left(B^\alg_{s-\delta}\right) \, 0_E, \quad  \e=
  \frac{\delta^{2k+2}}{4^{3k+4} c N(j)^2}, $$ 
 pour tout $\dt \in ]0,s[$ tel que $\dt \leq 4 N(j)^{1/(k+1)}$.}

\parag{5} L'hypoth{\`e}se sur $j$ et la condition $(A_c)$ permettent de
d{\'e}finir une application 
$$\phi : \bigcup_{s,\sigma} \left(\sigma B^{\alg}_{s+2\sigma}\right)
\rightarrow \alg, \quad \xi \mapsto j \, \left(e^{\xi}\right)^{-1} \,
\xi \, 0_E$$ 
et montrent que l'on a~:
$$|\phi(\xi)|_s \leq c N(j) \sigma^{-k-1} |\xi|_{s+2\sigma}^2 \quad
$$
pour tout $ \xi \in \sigma B^\alg_{s+2\sigma}$.  Fixons $s<1$ et
$\delta <s$ avec $\delta < 4 N(j)^{1/(k+1)}$. Soit $x \in \e
B^E_{s}$, $\e>0$, nous allons montrer que si $\e$ est choisit
suffisamment  petit alors les images successives de $j(x)$ par $\phi$,
soit 
$$\xi_n := \phi^n(j(x)),$$
sont bien d{\'e}finies pour tout $n\geq0$ et tendent vers $0$ dans
$B^\alg_{s-\delta}$. Nous montrerons ensuite que la suite $(e^{\xi_0}
\, ... \, e^{\xi_n})_{n\geq 0}$ tend vers un {\'e}l{\'e}ment $g$ dans $\exp
(B^\alg_{s-\delta})$. La suite d{\'e}finie par~:
$$x_0=x, \quad x_{n+1} = \left(e^{\xi_n}\right)^{-1} \, x_n,$$
tend alors vers $0_E$ dans $E_{(s-\delta)/2}$ de sorte que $x = g \,
0_E$. 

Soit donc $(s_n)$ la suite d{\'e}croissant vers $s-\delta$ d{\'e}finie par~:
$$s_0=s, \quad s_{n+1}=s_n-2\sigma_n \ {\rm avec\ } \sigma_n = 2^{-(n+2)}
\delta,$$
et posons~:
$$\e = \frac{c \,\delta}{16 \sigma_0} \prod_{m\geq 0} (c N(j)
\sigma_m^{-k-1})^{-2^{-m}}.$$ 
Nous v{\'e}rifierons, par la suite, que cette valeur de $\e$ co{\"\i}ncide avec
celle de l'{\'e}nonc{\'e}. 

Remarquons tout d'abord l'on a
$\left(2^{-4}\delta\right)^{2^n}\leq \sigma_{n+1}$ pour tout $n \geq 0$.  Pour
$n=0$ c'est {\'e}vident et une r{\'e}currence sur $n$
donne $$(2^{-4}\delta)^{2^{n}} =
\left((2^{-4}\delta)^{2^{n-1}}\right)^2 \leq (\sigma_n)^2 = 
2^{-2(n+2)} \delta \leq \sigma_{n+1}.$$ 

Supposons construits les termes $\xi_0,\dots,\xi_{N}$ de la suite.  Le
lemme suivant montre que $\xi_N$ appartient au domaine de d{\'e}finition
de $\phi$ et que l'on peut par cons{\'e}quent construire la suite
$(\xi_n)$ de proche en proche.

\begin{lemme}
  \label{lm:xi}
  Pour tout $0\leq n\leq N$, on a l'in{\'e}galit{\'e}~: 
  $$|\xi_n|_{s_{n+1}}\leq \left(2^{-4}\delta\right)^{2^n}.$$
\end{lemme}
\begin{proof}
  Nous allons montrer, par r{\'e}currence sur $n$, l'estimation~: 
  $$|\xi_n|_{s_{n+1}} \leq (2^{-4} \delta \, \mu_n)^{2^n} {\rm\ avec\ } \mu_n:=
  \prod_{m \geq n+1} (c N(j) \sigma_m^{-k-1})^{-2^{-m}}.$$ Cette
  estimation est plus pr{\'e}cise que celle du lemme car, d'apr{\`e}s
  l'hypoth{\`e}se sur $\delta$, on a $c N(j) \sigma_m^{-k-1} \geq 1$ pour
  tout $ m$, donc $\mu_n \leq 1$.  Pour $n=0$, l'hypoth{\`e}se sur $j$
  appliqu{\'e}e {\`a} $\xi_0=j(x)$ montre que l'on a~:
  $$|\xi_0|_{s_1} \leq N(j) \s_0^{-k} |x|_{s_1+\s_0} \leq N(j)
  \s_0^{-k} |x|_{s_0} \leq N(j) \sigma_0^{-k} \e = 2^{-4} \delta \,
  \mu_0.$$ 
  Supposons maintenant que
  $|\xi_{n-1}|_{s_n} \leq (2^{-4}\delta\, \mu_{n-1})^{2^{n-1}}$.\\
  L'hypoth{\`e}se sur $j$ et la conditionne $(A_c)$ donnent alors l'estimation~:
  $$\left| \xi_n \right|_{s_{n+1}} = \left| \phi(\xi_{n-1})
  \right|_{s_{n+1}} \leq c N(j) \s_n^{-k-1} |\xi_{n-1}|_{s_n}^2.$$
 En utilisant l'hypoth{\`e}se de r{\'e}currence, on obtient~:
  $$\left| \xi_n \right|_{s_{n+1}} \leq c N(j) \s_n^{-k-1}
  \left( 2^{-4} \delta \mu_{n-1} \right)^{2^n} = \left(2^{-4} \delta
    \, \mu_n \right)^{2^n},$$ 
  ce qui ach{\`e}ve la d{\'e}monstration du lemme.
\end{proof}

\parag{6} Nous allons montrer que la suite d{\'e}finie par~:
$$\gamma_0:= \xi_0, \quad \gamma_{n+1} := \log \left( e^{\gamma_{n}} \,
  e^{\xi_{n+1}} \right)$$
peut {\^e}tre construite de proche en proche et qu'elle converge dans la
boule $B^\alg_{s-\delta}$. 

Supposons construits les termes $\g_0,\dots,\g_N$ de la suite. Comme
le groupe $G$ est {\'e}chelonn{\'e}, la loi interne induit une application
$$B^\alg_{s_{n+1}} \times \sigma_{n+1} B^\alg_{s_{n+2}} \rightarrow
B^\alg_{s_{n+2}}, \quad (\gamma,\xi) \mapsto \log e^\gamma \, e^\xi.$$ 
Le lemme suivant montre que la suite $(\g_n)$ peut-{\^e}tre
construite de proche en proche dans la boule $B^\alg_{s-\delta}$.

\begin{lemme}
\label{lm:gamma}
  Pour tout $n\leq N$, on a l'in{\'e}galit{\'e}  $|\gamma_n|_{s_{n+1}} \leq 1$. 
\end{lemme}
\begin{proof}
  Nous allons montrer, par r{\'e}currence sur $n$, l'estimation
  $|\gamma_n|_{s_{n+1}} \leq g_n$, avec 
  $$g_0 := 2^{-4}, \quad g_{n} := \left(1 + 2^{-2^{n+1}}
  \right) g_{n-1} + 2^{-2^{n+1}}.$$ Cette estimation est plus pr{\'e}cise
  que celle du lemme parce que $(g_n)$ est major{\'e}e par $1$. En effet,
  une r{\'e}currence sur $n$ montre que l'on a~:
  $$1 + g_{n}=\prod_{1\leq m \leq n}
  \left(1+2^{-2^{m+1}}\right),$$ et par cons{\'e}quent~:
  $$\log (1+g_{n}) \leq 2^{-4} + \sum_{m\geq 1} 2^{-2^{m+1}}\leq
  2^{-4} +\sum_{m\geq 1} 2^{-4^m}= \frac{1}{16} + \frac{1}{28}<\log
  2.$$

  Puisque $\gamma_0=\xi_0$, l'estimation voulue est vraie au rang
  $n=0$ d'apr{\`e}s le lemme~\ref{lm:xi}. Supposons la propri{\'e}t{\'e}
  satisfaite au rang $n-1$. Commme $G$ est {\'e}chelonn{\'e}, on
  a~: $$|\gamma_n|_{s_{n+1}} \leq |\xi_n|_{s_{n+1}} + \sigma_n^{-1}
  |\xi_n|_{s_{n+1}} |\gamma_{n-1}|_{s_n} + |\gamma_{n-1}|_{s_n}.$$
  Dans le membre de droite, d'apr{\`e}s le lemme~\ref{lm:xi} et
  l'hypoth{\`e}se de r{\'e}currence, les trois termes peuvent {\^e}tre major{\'e}s
  respectivement par $2^{-2^{n+1}}$, $2^{-2^{n+1}} g_{n-1}$ et
  $g_{n-1}$, et nous obtenons donc~:
  $$|\gamma_n|_{s_{n+1}} \leq \left(1 + 2^{-2^{n+1}} \right) g_{n-1} +
  2^{-2^{n+1}} = g_{n},$$ ce qui d{\'e}montre le lemme.
\end{proof}

\begin{lemme}
  La suite $(\gamma_n)$ converge dans la boule $B^\alg_{s-\delta}$.
\end{lemme}

\begin{proof}
  En substituant $\log \left( e^{\gamma_{n-1}} \, e^{\xi_{n}} \right)$
  {\`a} $\g_n$ puis en appliquant l'in{\'e}galit{\'e} triangulaire, nous
  obtenons~: 
  $$|\gamma_n-\gamma_{n-1}|_{s_{n+1}} \leq |\xi_n|_{s_{n+1}} + |\log
  e^{\gamma_{n-1}} \, e^{\xi_n} - \gamma_{n-1} - \xi_n|_{s_{n+1}}.$$
  En utilisant le fait que $G$ est {\'e}chelonn{\'e}, le membre de droite se
  simplifie~:
  $$|\gamma_n-\gamma_{n-1}|_{s_{n+1}} \leq |\xi_n|_{s_{n+1}} + \kappa
  \sigma_n^{-1} |\gamma_{n-1}|_{s_n} |\xi_n|_{s_{n+1}},$$
  et par cons{\'e}quent, en utilisant le lemme~\ref{lm:xi}, il vient~: 
  $$|\gamma_n-\gamma_{n-1}|_{s -\dt}
  \leq|\gamma_n-\gamma_{n-1}|_{s_{n+1}} \leq 2^{-2^{n+2}} + \kappa
  2^{-2^{n+1}} \leq (1+\kappa) 2^{-2^{n+1}},$$ d'o{\`u} la convergence
  souhait{\'e}e.
\end{proof}
  
Il ne nous reste plus qu'{\`a} calculer $\e$. En substituant
$2^{-(n+2)}\delta$ {\`a} $\sigma_n$ dans la d{\'e}finition de $\e$ on
obtient~: 
$$\e = \frac{c}{4} \prod_{m\geq 0} \left(
  \frac{\delta^{k+1}}{cN(j)2^{(k+1)(m+2)}} \right)^{2^{-m}},$$
soit
$$\log \e = \log \frac{c}{4} + \left(\sum_{m\geq 0} 2^{-m} \right)
\log \frac{\delta^{k+1}}{cN(j)4^{k+1}} - \left( \sum_{m\geq 0} m
  2^{-m} \right) \log 2^{k+1}$$
ce qui donne bien la valeur annonc{\'e} car $\sum_{m\geq 0} 2^{-m} =
\sum_{m\geq 0} m2^{-m} =2$. Le th{\'e}or{\`e}me est d{\'e}montr{\'e}.

\end{document}